\pdfoutput=1
\documentclass[reqno,12pt]{amsart}
\usepackage[T1]{fontenc}
\usepackage[utf8]{inputenc}
\usepackage{lmodern}
\usepackage[english]{babel}
\usepackage[margin=25mm]{geometry}
\usepackage{amsmath,amssymb,amsthm,mathtools}
\usepackage{booktabs,array,tabularx}
\usepackage{enumitem}
\usepackage{microtype}
\usepackage{tikz}

\usetikzlibrary{positioning,fit,backgrounds,calc}
\newtheorem{theorem}{Theorem}
\newtheorem{lemma}{Lemma}

\newcommand{\partner}[1]{\overline{#1}}
\newcommand{\mustar}{\mu^*}
\usepackage{needspace}
\usepackage[hidelinks]{hyperref}

\allowdisplaybreaks[1]

\newcommand{\bricks}{\mathcal J}
\newcommand{\anch}{\mathcal A}

\usepackage{etoolbox}
\makeatletter
\patchcmd{\@setauthors}{\@topsep30\p@}{\@topsep15\p@}{}{}
\patchcmd{\@setauthors}{\MakeUppercase{\authors}}{\authors}{}{}
\makeatother

\begin{document}

\title[The saturation number is at most the harmonic index]{The saturation number is at most \\the harmonic index
in subquartic graphs}
\author[]{%
  {\normalfont\fontsize{8.0}{10.0}\selectfont Jakub Zabrzeski}\\[-1pt]
  {\normalfont\fontsize{6.8}{8.8}\selectfont Faculty of Mathematics and Computer Science, Adam Mickiewicz University, Poznań, Poland}%
}
\date{27 November 2026}

\begin{abstract}
For a finite simple graph $G$, let $\mustar(G)$ denote the minimum
cardinality of a maximal matching and let
$H(G)=\sum_{uv\in E_G}2/(d(u)+d(v))$ be its harmonic index.
We prove that $\mustar(G)\le H(G)$ whenever $\Delta(G)\le4$,
resolving Gupta's conjecture for subquartic graphs. 
\end{abstract}


\vspace*{-15pt}
\maketitle

\section{Introduction}
A matching $M$ in a graph $G=(V_G,E_G)$ is {\bf maximal} if no edge can be added to it.
The minimum cardinality of such a matching is the {\bf saturation number} of $G$, denoted as 
$\mustar(G)$. The harmonic index of $G$ is defined as 
\[
 H(G)=\sum_{uv\in E_G}\frac{2}{d(u)+d(v)}\quad \textrm{and}\quad H(G)=0\textrm{\ for\ }E_G=\emptyset\,. 
\]
Thus, for instance, $\mustar(K_2)=H(K_2)=1$, while for the path $P_3$ on three vertices 
$$\mustar(P_3)=1\le H(P_3)=4/3\,.$$
Gupta~\cite{Gupta} asked for which graphs $G$ we have 
$\mustar(G)\le H(G)$, and speculated (\cite[Section 5, Conjecture 1]{Gupta}) that it holds for each graph 
$G$ with $\Delta(G)\le 4$. Our main result settles his conjecture in the affirmative.

\begin{theorem}\label{thm:main}
If $\Delta(G)\le4$, then $\mustar(G)\le H(G)$.
\end{theorem}

Let us note that this result is, in a way, best possible. Indeed, if 
$S^*_5$ denote the `subdivided star' with vertex set
$\{c\}\cup\{u_i,\ell_i:1\le i\le 5\}$ and edge set
$\{cu_i,u_i\ell_i:1\le i\le 5\}$; then
\[
 \Delta(S^*_5)=5,\qquad \mustar(S^*_5)=5,\qquad H(S^*_5)=\frac{100}{21}<5.
\]
Thus, in Theorem~\ref{thm:main} the bound for $\Delta(G)$ is best possible and 
$\Delta(G)$ cannot be replaced by the average degree. 

The proof follows a rather natural idea although its details are somewhat technical:
initially we assign a unit weight to each edge of the maximal matching and then distribute them among all edges of the graph so that new weights of the edges $e=vw$ can be easily compared to quotient $2/(d(v)+d(w))$ 
which appear in the definition of the harmonic number.
Thus, in Section 2 we prove some structural results on graphs with maximal matchings of a given size. 
Next, we construct the vertex labels used to define the auxiliary vertex weights. 
Section~\ref{sec:weights} derives the weight functions from the local estimates they must satisfy. Section~5 shows how matching edges near unsaturated vertices cover
their cost while 
Section~6 verifies that edges joining matching edges supply
the required amounts assigned to their endpoints.
Section~7 handles the remaining edges. Finally, in the Section~8 we collect all above estimates completing the proof of Theorem~\ref{thm:main}.

\section{The matching and its induced components}
In this section we study the structure of graphs with a given saturation number.
For a matching $M$ by $U_M=V_G\setminus V_M$ we denote 
the set of vertices unsaturated by $M$. Note that $M$ is maximal if and only if 
$U_M$ is independent -- we shall often use this observation throughout this note. 
In our proof we choose $M$ among all maximal matchings of minimum
cardinality so as to minimize
\begin{equation}\label{eq00}
 \sum_{z\in U_M}d(z).
\end{equation}
For a saturated vertex $v$, we denote its unique
neighbour in $M$ by $\partner v$, and put
\[
 t_v=|N(v)\cap U_M|.
\]
We call a saturated vertex {\bf clean} if $t_v=0$ and {\bf active} otherwise.

\begin{lemma}[Alternating-path exchange]\label{lem:path}
If an $M$-alternating path of even length, starts at $z\in U_M$
and ends at a clean vertex $v$, then $d(v)\ge d(z)$.
\end{lemma}
\begin{proof}
Consider an $M$-alternating path $x_0x_1\ldots x_{2k}$, where $x_0=z$,
$x_{2k}=v$, and $k\ge1$. Let
\[
 E^-:=\{x_{2i-1}x_{2i}:1\le i\le k\}\subseteq M,\qquad
 E^+:=\{x_{2i-2}x_{2i-1}:1\le i\le k\},
\]
and $M'=(M\setminus E^-)\cup E^+$. Then $|M'|=|M|$ and 
\[
 V_M\setminus V_{M'}=\{v\},\quad
 V_{M'}\setminus V_M=\{z\},\quad U_{M'}=(U_M\setminus\{z\})\cup\{v\}.
\]
Since the set $U_M\setminus\{z\}$ is independent and $v$ has no neighbor
in it because $v$ is clean with respect to $M$ the set $U_{M'}$ is independent, and so $M'$ is maximal. Moreover, by our choice of~$M$, 
\[
 \sum_{u\in U_M}d(u)\le\sum_{u\in U_{M'}}d(u)
 =\sum_{u\in U_M}d(u)-d(z)+d(v).
\]
Hence $d(v)\ge d(z)$, as required.
\end{proof}

Partition $M$ into $P$ and $Q$: an edge belongs to $P$ if at least
one endpoint is active, and to $Q$ otherwise. Every edge
incident with $U_M$ has its other endpoint in $V_P$. Next we determine the components of $G[V_Q]$, and prove two exchange lemmas about edges entering $V_Q$.

\Needspace{8\baselineskip}

\begin{lemma}\label{lem:bricks}
Each component of $G[V_Q]$ is isomorphic to $K_2$, $C_4$, $K_4$, $K_{3,3}$ or $K_{4,4}$.
\end{lemma}

\begin{proof}
Let $J$ be a component of $G[V_Q]$; then $Q\cap E(J)$ is a perfect matching of $J$. If $J$ has a maximal matching $M_J$ of size less than $|V_J|/2$, replacing $Q\cap E(J)$ by $M_J$ in $M$ leaves $U_M\cup(V_J\setminus V(M_J))$ unsaturated, which is independent (no edge joins $V_Q$ to $U_M$) and contradicts the minimality of $|M|$. Hence every maximal matching of $J$ is perfect. From Sumner's theorem~\cite{Sumner}, $J$ is $K_{2k}$ or $K_{k,k}$ ($k\ge 1$), and $\Delta(G)\le 4$ leaves only the five graphs listed. This proves Lemma~\ref{lem:bricks}.
\end{proof}

We call the components of $G[V_Q]$ \emph{$Q$-bricks} (or \emph{bricks}), and let $\bricks$ be their family. Each $J\in\bricks$ is internally $r_J$-regular ($1\le r_J\le 4$), and every edge of $J$ lies in a perfect matching of $J$. Replacing the matching on a brick by any other perfect matching of it leaves $|M|$, $U_M$, the matching on $V_P$, and the sum in~\eqref{eq00} unchanged. We use this freedom in the next two lemmas.

\Needspace{7\baselineskip}

\begin{lemma}\label{lem:two-to-one}
Let $ab, xy$ be distinct edges of $M$. If $ax\in E_G$ and both $b, y$ are clean, then $by\in E_G$.
\end{lemma}

\begin{proof}
Suppose not, and let $M'=(M\setminus\{ab,xy\})\cup\{ax\}$. The set of vertices unsaturated by $M'$ is $U_M\cup\{b,y\}$. Since $U_M$ is independent, both $b, y$ are clean, and $by\notin E_G$, it follows that $U_M\cup\{b,y\}$ is independent; so $M'$ is a maximal matching of $G$ smaller than $M$, a contradiction. This proves Lemma~\ref{lem:two-to-one}.
\end{proof}

It follows immediately that if $ab\in P$ has $b$ clean, $J\in\bricks$ is a brick, $x\in V_J$, and $ax\in E_G$, then
\begin{equation}\label{eq:brick-nbr}
N_J(x)\subseteq N(b).
\end{equation}
To see this, let $y\in N_J(x)$, replace the matching on $J$ by a perfect matching of $J$ containing $xy$, and apply Lemma~\ref{lem:two-to-one}; we deduce that $by\in E_G$. Next, we show that an edge of $P$ cannot attach to two different bricks from opposite ends.

\begin{lemma}\label{lem:same-brick}
If both ends of an edge $ab\in P$ have neighbours in $V_Q$, then all their neighbours in $V_Q$ lie in the same brick.
\end{lemma}

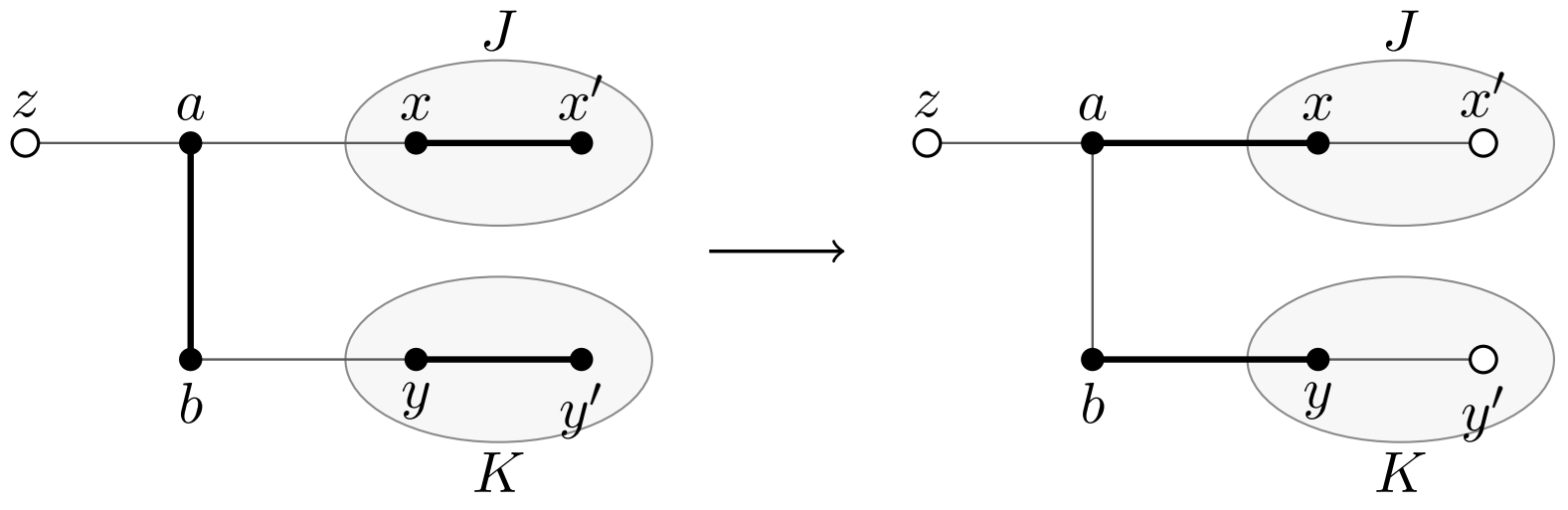
\begin{figure}[ht]
\centering
\begin{tikzpicture}[x=1cm,y=1cm,font=\small,
 dot/.style={circle,draw=black,fill=black,inner sep=0pt,minimum size=4pt},
 free dot/.style={dot,fill=white,minimum size=5pt,line width=.6pt},
 thin edge/.style={draw=black!65,line width=.45pt},
 matching edge/.style={draw=black,line width=1.2pt},
 brick outline/.style={draw=black!45,fill=black!3,line width=.4pt},
 every label/.style={inner sep=2pt}]
\begin{scope}
 \path[brick outline] (3.15,.72) ellipse (1.02cm and .55cm);
 \path[brick outline] (3.15,-.72) ellipse (1.02cm and .55cm);
 \node at (3.15,1.47) {$J$};
 \node at (3.15,-1.47) {$K$};
 \coordinate (z) at (0,.72);
 \coordinate (a) at (1.1,.72);
 \coordinate (b) at (1.1,-.72);
 \coordinate (x) at (2.6,.72);
 \coordinate (xp) at (3.7,.72);
 \coordinate (y) at (2.6,-.72);
 \coordinate (yp) at (3.7,-.72);
 \draw[thin edge] (z)--(a) (a)--(x) (b)--(y);
 \draw[matching edge] (a)--(b) (x)--(xp) (y)--(yp);
 \node[free dot,label=above:$z$] at (z) {};
 \node[dot,label=above:$a$] at (a) {};
 \node[dot,label=below:$b$] at (b) {};
 \node[dot,label=above:$x$] at (x) {};
 \node[dot,label=above:$x'$] at (xp) {};
 \node[dot,label=below:$y$] at (y) {};
 \node[dot,label=below:$y'$] at (yp) {};
 \node at (2,-2.0) {(a) $M$};
\end{scope}
\draw[->,line width=.65pt] (4.55,0)--(5.45,0);
\begin{scope}[xshift=6cm]
 \path[brick outline] (3.15,.72) ellipse (1.02cm and .55cm);
 \path[brick outline] (3.15,-.72) ellipse (1.02cm and .55cm);
 \node at (3.15,1.47) {$J$};
 \node at (3.15,-1.47) {$K$};
 \coordinate (z) at (0,.72);
 \coordinate (a) at (1.1,.72);
 \coordinate (b) at (1.1,-.72);
 \coordinate (x) at (2.6,.72);
 \coordinate (xp) at (3.7,.72);
 \coordinate (y) at (2.6,-.72);
 \coordinate (yp) at (3.7,-.72);
 \draw[thin edge] (z)--(a) (a)--(b) (x)--(xp) (y)--(yp);
 \draw[matching edge] (a)--(x) (b)--(y);
 \node[free dot,label=above:$z$] at (z) {};
 \node[dot,label=above:$a$] at (a) {};
 \node[dot,label=below:$b$] at (b) {};
 \node[dot,label=above:$x$] at (x) {};
 \node[free dot,label=above:$x'$] at (xp) {};
 \node[dot,label=below:$y$] at (y) {};
 \node[free dot,label=below:$y'$] at (yp) {};
 \node at (2,-2.0) {(b) $M'$};
\end{scope}
\end{tikzpicture}
\caption{The three-to-two exchange
$M'=(M\setminus\{ab,xx',yy'\})\cup\{ax,by\}$ in Lemma~\ref{lem:same-brick}.
Thick edges belong to the displayed matching; open circles mark
unsaturated vertices. Only the relevant vertices and edges of the distinct
bricks $J,K$ are shown. The roles of $a,b$ may be interchanged so that
$a$ has an unsaturated neighbour $z$.}
\label{fig:exchange}
\end{figure}

\begin{proof}
Suppose that $ax, by\in E_G$, where $x\in V_J$ and $y\in V_K$ for distinct bricks $J, K\in\bricks$. Choose internal neighbours $x'\in N_J(x)$ and $y'\in N_K(y)$, and replace the matchings on $J, K$ by perfect matchings containing $xx', yy'$ respectively. Let $M'=(M\setminus\{ab,xx',yy'\})\cup\{ax,by\}$ (see Figure~\ref{fig:exchange}). Then $|M'|=|M|-1$, and the vertices unsaturated by $M'$ are $U_M\cup\{x',y'\}$. Since $x', y'$ are clean and belong to distinct components of $G[V_Q]$, they are nonadjacent to each other and to every vertex of $U_M$. Hence $U_M\cup\{x',y'\}$ is independent, so $M'$ is a maximal matching of $G$ smaller than $M$, a contradiction.

Now fix a neighbour of $a$ in $V_Q$, in some brick $J$ say. By what we have just proved, every neighbour of $b$ in $V_Q$ also lies in $J$; and fixing one such neighbour of $b$, the same argument shows that every neighbour of $a$ in $V_Q$ lies in $J$ as well. This proves Lemma~\ref{lem:same-brick}.
\end{proof}

Let us say that an edge $ab\in P$ is \emph{anchored} if both $a, b$ have neighbours in $V_Q$. By Lemma~\ref{lem:same-brick}, all these neighbours belong to a unique brick, to which we assign $ab$.

\section{Labels propagated by alternating paths}\label{sec:labels}

We need to keep track of the degree bounds coming from Lemma~\ref{lem:path}, and it is convenient to record them as vertex labels. For each saturated vertex $v$, put
\[
h_v = \max\{d(z) : z\in N(v)\cap U_M\} \quad \text{if } t_v>0,
\]
and $h_v=0$ if $t_v=0$. Let $T$ be the set of clean ends of edges in $P$; so each $v\in T$ has an active partner $\partner{v}$. We define a digraph $\mathcal{D}$ with vertex set $T$ by putting an arc $v\to u$ whenever $v, u\in T$ are distinct and $v\partner{u}\in E_G$.

We assign labels to the vertices of $T$ as follows: each $v\in T$ receives the \emph{initial label} $h_{\partner{v}}$, and whenever $ab\in P$ has both ends active, we also give the \emph{additional label} $h_b$ to every $v\in T$ with $a\partner{v}\in E_G$ (and symmetrically $h_a$ to every $v\in T$ with $b\partner{v}\in E_G$). For $v\in T$, let $R_v$ be the largest label that reaches $v$ along a directed path of $\mathcal{D}$ (including paths of length zero).

\begin{lemma}\label{lem:labels}
If $v\in T$, then $1\le h_{\partner{v}}\le R_v\le d(v)\le 4$.
\end{lemma}

\begin{proof}
Since $\partner{v}$ is active, it has a neighbour in $U_M$, so $h_{\partner{v}}\ge 1$; and the zero-length path at $v$ gives $R_v\ge h_{\partner{v}}$. Now let $v_0,v_1,\dots,v_k=v$ be a directed path in $\mathcal{D}$ (with $v_0,\dots,v_k$ all distinct) carrying a label $d(z)$ from $v_0$ to $v$. Either this is the initial label $h_{\partner{v}_0}=d(z)$ with $z\in N(\partner{v}_0)\cap U_M$, in which case we take the $M$-alternating path
\[
z,\partner{v}_0,v_0,\partner{v}_1,v_1,\dots,\partner{v}_k,v_k,
\]
or it is an additional label $h_b=d(z)$ from an edge $ab\in P$ with both ends active, $a\partner{v}_0\in E_G$, and $z\in N(b)\cap U_M$, in which case we take
\[
z,b,a,\partner{v}_0,v_0,\partner{v}_1,v_1,\dots,\partner{v}_k,v_k.
\]
Here $v_0,\dots,v_k$ are distinct clean vertices, their partners $\partner{v}_0,\dots,\partner{v}_k$ are distinct active vertices, and $z\in U_M$; and in the second case neither $a$ nor $b$ belongs to $\{v_0,\dots,v_k,\partner{v}_0,\dots,\partner{v}_k\}$ because both $a, b$ have active partners. Thus this path is simple, of even length, and ends at the clean vertex $v$. By Lemma~\ref{lem:path}, $d(v)\ge d(z)$. Taking the maximum over all labels reaching $v$ gives $R_v\le d(v)\le 4$. This proves Lemma~\ref{lem:labels}.
\end{proof}

For the weight estimates that follow, it is convenient to attach a single number to each vertex $v\in V_P$. Let us define the \emph{signal} $s_v$ of $v\in V_P$ by
\begin{equation}\label{eq:signal}
s_v = \begin{cases}
R_v & \text{if } v \text{ is clean,}\\
h_{\partner{v}} & \text{if } v \text{ is active.}
\end{cases}
\end{equation}
Thus, $s_v=0$ if and only if $v$ is active and $\partner{v}$ is clean (that is, $v$ has unsaturated neighbours and $\partner{v}$ does not); and $s_v>0$ for every other vertex $v\in V_P$. (Later we will have to treat the zero-signal vertices separately, with a different weight function.) First, we need the following bounds on degrees across edges between $P$-edges, and across boundary edges into bricks.

\Needspace{23\baselineskip}

\begin{lemma}\label{lem:signal}
\begin{enumerate}
\item[\rm(a)] If $a, c\in V_P$ belong to distinct edges of $P$, and $a\partner{c}\in E_G$, and $c\in T$, then $s_a\le R_c$.
\item[\rm(b)] If $ax\in E_G$, where $a\in V_P$ and $x\in V_J$ for some brick $J\in\bricks$, then
\begin{equation}\label{eq:brick-signal}
s_a\le \min_{y\in N_J(x)} d(y).
\end{equation}
\end{enumerate}
\end{lemma}

\begin{proof}
If $s_a=0$ then both statements hold trivially (since $R_c\ge 1$ by Lemma~\ref{lem:labels}, and every $y\in N_J(x)$ has degree at least one). So we may assume that $s_a>0$; hence either $a\in T$, or both ends of $a\partner{a}$ are active.

For~(a), if $a\in T$, then $a\to c$ is an arc of $\mathcal{D}$, so every label reaching $a$ also reaches $c$, and therefore $s_a=R_a\le R_c$. If both ends of $a\partner{a}$ are active, then $s_a=h_{\partner{a}}$ is one of the additional labels placed at $c$, and again $s_a\le R_c$.

For~(b), let $y\in N_J(x)$, and replace the matching on $J$ by a perfect matching of $J$ containing $xy$. If $a\in T$, take an $M$-alternating path ending at $a$ that realises $R_a$, as in the proof of Lemma~\ref{lem:labels}, and extend it by adding the edges $ax, xy$; all its earlier saturated vertices are in $V_P$, so this path is simple. If both ends of $a\partner{a}$ are active, choose $z\in N(\partner{a})\cap U_M$ with $d(z)=h_{\partner{a}}=s_a$, and take the simple $M$-alternating path $z,\partner{a},a,x,y$. In either case $y$ is clean, so $d(y)\ge s_a$ by Lemma~\ref{lem:path}. Since this holds for every $y\in N_J(x)$,~\eqref{eq:brick-signal} follows. This proves Lemma~\ref{lem:signal}.
\end{proof}

\section{Degree weights}\label{sec:weights}

Now we define the local weights. First we look at a single brick: we know the cost of its perfect matching, and we want to see how much its boundary edges must supply to make up for the deficits of its internal edges. Whatever weight is left over on a boundary edge then tells us how much we can ask for at its other end, and this gives us our ordinary weight function $f$. (When the signal is zero, we will need a separate weight function $g$, defined later in this section.)

Let $J\in\bricks$ be a brick, with $n=|V_J|$ vertices and internal degree $r=r_J$ say. The matching $Q\cap E(J)$ has $n/2$ edges, whereas $J$ has $rn/2$ internal edges, so dividing the cost equally gives each internal edge a target weight of $1/r$. At an internal edge $xy\in E_J$, the deficit is
\[
\frac{1}{r}-\frac{2}{d(x)+d(y)}=\frac{d(x)-r}{r(d(x)+d(y))}+\frac{d(y)-r}{r(d(x)+d(y))}.
\]
We charge the first term to $x$ and the second to $y$. If $d(x)>r$, then $x$ has $d(x)-r$ boundary edges; dividing the charge at $x$ equally among them asks each boundary edge at $x$ for $1/[r(d(x)+d(y))]$. Summing over all $y\in N_J(x)$, each boundary edge incident with $x$ is asked to supply
\[
\phi_{J,x}=\frac{1}{r}\sum_{y\in N_J(x)}\frac{1}{d(x)+d(y)}.
\]
(Here $r$ is the internal degree of $J$, while $d(x), d(y)$ are full degrees in $G$.) The total contribution collected at $x$ is $(d(x)-r)\phi_{J,x}$ (which is zero if $d(x)=r$); and summing the edge deficits over all $xy\in E_J$ gives
\begin{equation}\label{eq:brick-balance}
\sum_{xy\in E_J}\frac{2}{d(x)+d(y)}-\frac{|V_J|}{2}+\sum_{x\in V_J}(d(x)-r)\phi_{J,x}=0.
\end{equation}

Now let $ax$ be a boundary edge, where $a\in V_P$, $x\in V_J$, and $s_a>0$. By Lemma~\ref{lem:signal}(b), each of the $r$ summands in $\phi_{J,x}$ is at most $1/(d(x)+s_a)$, and so
\[
\phi_{J,x}\le \frac{1}{r}\cdot \frac{r}{d(x)+s_a}=\frac{1}{d(x)+s_a}.
\]
After paying its share to the brick, the edge $ax$ still has at least
\[
\frac{2}{d(a)+d(x)}-\phi_{J,x}\ge \frac{2}{d(a)+d(x)}-\frac{1}{d(x)+s_a}
\]
left over for $a$. Let $d=d(a)$, $s=s_a$, and $k=d(x)$. Since $x$ has an internal neighbour in $J$ and the boundary neighbour $a$, we have $2\le k\le 4$. To get a weight depending only on $d, s$ that works for every possible degree $k$, we take the minimum: for integers $2\le d\le 4$ and $1\le s\le 4$, define
\begin{equation}\label{eq:f-def}
f(d,s)=\min_{2\le k\le 4}\left(\frac{2}{d+k}-\frac{1}{k+s}\right).
\end{equation}
We can evaluate this explicitly:
\begin{equation}\label{eq:f-explicit}
f(d,s)=\begin{cases}
(4-d)/15 & \text{if } s=1,\\[2pt]
2/(d+4)-1/(s+4) & \text{if } 2\le s\le 4.
\end{cases}
\end{equation}
To see this, when $s=1$, checking $k\in\{2,3,4\}$ directly shows that the minimum occurs at $k=4$ if $d=2$, and at $k=2$ if $d\in\{3,4\}$. When $s\ge 2$, subtracting the value at $k=4$ from the value at $k$ gives
\[
(4-k)\left(\frac{2}{(d+k)(d+4)}-\frac{1}{(s+k)(s+4)}\right).
\]
Since $d\le 4$ and $s\ge 2$, the expression in parentheses is at least
\[
\frac{1}{4(k+4)}-\frac{1}{6(k+2)}=\frac{k-2}{12(k+4)(k+2)}\ge 0,
\]
so the minimum is attained at $k=4$, and~\eqref{eq:f-explicit} follows. From~\eqref{eq:f-def}, $f(d,s)$ is nondecreasing in $s$, and therefore
\begin{equation}\label{eq:u-bound}
0\le f(d,s)\le u_d := f(d,4)=\frac{2}{d+4}-\frac{1}{8}.
\end{equation}

What if an end $a\in V_P$ has signal $s_a=0$? Then $a$ is active and its partner $b=\partner{a}$ is clean, so we have the propagated label $R=R_b\in\{1,2,3,4\}$ at $b$. For this case we define the \emph{exceptional weights} by
\begin{equation}\label{eq:g3}
g(3,R)=\frac{5}{7}-\frac{2}{R+3}-3f(4,R),\qquad 1\le R\le 4,
\end{equation}
\begin{equation}\label{eq:g4}
g(4,1)=\frac{11}{60},\qquad g(4,R)=\frac{1}{2(R+4)},\qquad 2\le R\le 4.
\end{equation}
(We will see in Sections~\ref{sec:matchingcost} and~\ref{sec:anchors} where these values come from: they make the tightest local inequalities hold with equality.) All values of $f, g$ are given in Table~\ref{tab:weights}. Checking the table shows that $0\le g(d,R)<2/(d+R)$ in all cases, and $g(d,R)\le u_d$ whenever $R\ge 2$. When terms appear with coefficient zero, we put $g(2,R)=0$ and $f(1,1)=0$.

\begin{table}[htbp]
\centering
\renewcommand{\arraystretch}{1.2}
\begin{tabular}{ccccc}
\toprule
$s$ or $R$ & $1$ & $2$ & $3$ & $4$ \\
\midrule
$f(2,s)$ & $2/15$ & $1/6$ & $4/21$ & $5/24$ \\
$f(3,s)$ & $1/15$ & $5/42$ & $1/7$ & $9/56$ \\
$f(4,s)$ & $0$ & $1/12$ & $3/28$ & $1/8$ \\
$g(3,R)$ & $3/14$ & $9/140$ & $5/84$ & $3/56$ \\
$g(4,R)$ & $11/60$ & $1/12$ & $1/14$ & $1/16$ \\
\bottomrule
\end{tabular}
\caption{The ordinary and exceptional weights.}
\label{tab:weights}
\end{table}

Now we assign an auxiliary weight $\lambda_a$ to each end $a$ of an edge $ab\in P$. Let $c_a=d(a)-1-t_a$ be the number of saturated neighbours of $a$ other than $b$. If $c_a=0$, set $\lambda_a=0$; otherwise set
\[
\lambda_a=\begin{cases}
f(d(a),s_a) & \text{if } s_a>0,\\[2pt]
g(d(a),R_b) & \text{if } s_a=0.
\end{cases}
\]
We say that $a$ is \emph{exceptional} in the second case (that is, when $s_a=0$ and $c_a>0$). Note that if $a$ is exceptional, then $a$ is active, its partner $b$ is clean, and $d(a)=1+t_a+c_a\ge 3$, so $d(a)\in\{3,4\}$.

If $ax$ is an edge between a nonexceptional end $a\in V_P$ and a vertex $x\in V_J$ of some brick $J\in\bricks$, then $c_a>0$, and hence $s_a>0$; thus $\lambda_a=f(d(a),s_a)$, and from~\eqref{eq:f-def} we have
\begin{equation}\label{eq:ordinary-bdry}
\frac{2}{d(a)+d(x)}-\lambda_a-\phi_{J,x}\ge 0.
\end{equation}
So every boundary edge at a nonexceptional end covers both its auxiliary weight and its share to the brick. That leaves three things to check: the cost of matching edges in $P$, the edges between distinct matching edges of $P$, and the bricks attached to exceptional ends. We handle them in the next three sections.

\section{The cost of a matching edge}\label{sec:matchingcost}

First, we show that each edge $ab\in P$ covers its unit cost once it receives the auxiliary weights $\lambda_a, \lambda_b$ from its edges to other saturated vertices, on top of its own weight and the weights of its edges to $U_M$. For $ab\in P$, let
\[
\rho_{ab}=\frac{2}{d(a)+d(b)}-1+\sum_{\substack{v\in\{a,b\}\\t_v>0}}\frac{2t_v}{d(v)+h_v}.
\]
Here the sum is a lower bound on the weight of the edges between $\{a,b\}$ and $U_M$, and subtracting $1$ accounts for the matching edge $ab$ itself. We claim:

\begin{lemma}\label{lem:matching-cost}
For every edge $ab\in P$,
\begin{equation}\label{eq:Sab}
S_{ab}:=\rho_{ab}+c_a\lambda_a+c_b\lambda_b\ge 0.
\end{equation}
\end{lemma}

\begin{proof}
There are two cases, depending on whether both ends of $ab$ are active, or only one of them is.

Suppose first that both $a, b$ are active, with degrees $d=d(a)$ and $e=d(b)$ say. Then $s_a=h_b\ge 1$ and $s_b=h_a\ge 1$, so neither end is exceptional; thus $c_a\lambda_a=(d-1-t_a)f(d,h_b)$ and $c_b\lambda_b=(e-1-t_b)f(e,h_a)$. If we increase $t_a$ by one (leaving $d, e, h_a, h_b$ unchanged), then $S_{ab}$ increases by
\[
\frac{2}{d+h_a}-f(d,h_b)\ge \frac{2}{d+4}-u_d=\frac{1}{8}>0,
\]
and similarly for $t_b$. Consequently, $S_{ab}$ is at least its value when $t_a=t_b=1$, namely
\[
\frac{2}{d+e}-1+\left(\frac{2}{d+h_a}+(e-2)f(e,h_a)\right)+\left(\frac{2}{e+h_b}+(d-2)f(d,h_b)\right).
\]
For $2\le h\le 4$, we have $f(e,h)=2/(e+4)-1/(h+4)$, so as a function of $h$ the expression $2/(d+h)+(e-2)f(e,h)$ has derivative
\[
-\frac{2}{(d+h)^2}+\frac{e-2}{(h+4)^2}\le 0,
\]
because $d\le 4$ and $e-2\le 2$. Hence the minimum of $2/(d+h)+(e-2)f(e,h)$ over $h\in\{1,2,3,4\}$ is attained at $h=1$ or $h=4$. Evaluating the six unordered pairs $(d,e)$ gives the positive lower bounds in Table~\ref{tab:two-active}, so $S_{ab}>0$ in this case.

\begin{table}[ht]
\centering
\small
\setlength{\tabcolsep}{4pt}
\renewcommand{\arraystretch}{1.2}
\begin{tabular}{@{}cccc@{}}
\toprule
$(d,e)$ &
$\displaystyle\min_{h\in\{1,4\}}\left(\frac2{d+h}+(e-2)f(e,h)\right)$ &
$\displaystyle\min_{h\in\{1,4\}}\left(\frac2{e+h}+(d-2)f(d,h)\right)$ &
\shortstack{lower bound\\for $S_{ab}$}\\[5pt]
\midrule
$(2,2)$ & $1/3$ & $1/3$ & $1/6$ \\
$(2,3)$ & $83/168$ & $2/7$ & $151/840$ \\
$(2,4)$ & $7/12$ & $1/4$ & $1/6$ \\
$(3,3)$ & $25/56$ & $25/56$ & $19/84$ \\
$(3,4)$ & $1/2$ & $23/56$ & $11/56$ \\
$(4,4)$ & $2/5$ & $2/5$ & $1/20$ \\
\bottomrule
\end{tabular}
\caption{Lower bounds for $S_{ab}$ when both ends are active.}
\label{tab:two-active}
\end{table}

Now suppose that only one end of $ab$ is active: say $a$ is active of degree $d$, and $b$ is clean of degree $D$. Set $R=R_b$; by Lemma~\ref{lem:labels}, $1\le h_a\le R\le D\le 4$. Replacing $h_a$ by $R$ does not increase $2t_a/(d+h_a)$; and since $2/(d+R)>g(d,R)$, reducing $t_a$ to one cannot increase the bound either (noting that $c_a\lambda_a=(d-1-t_a)g(d,R)$, with our convention $g(2,R)=0$). It is therefore enough to prove that
\begin{equation}\label{eq:one-active}
\frac{2}{d+D}+\frac{2}{d+R}-1+(d-2)g(d,R)+(D-1)f(D,R)\ge 0.
\end{equation}
If $D=1$, then $R=1$, and the left-hand side of~\eqref{eq:one-active} becomes $4/(d+1)-1+(d-2)g(d,1)$, which equals $1/3$, $3/14$, or $1/6$ for $d=2,3,4$ respectively. So we may assume that $2\le D\le 4$.

For fixed $d, R$, we reduce to $D=4$. Subtracting the value of $2/(d+D)+(D-1)f(D,R)$ at $D=4$ from its value at $D$ gives $4-D$ times
\[
\begin{cases}
\displaystyle\frac{2}{(d+D)(d+4)}+\frac{D-1}{15}, & R=1,\\[6pt]
\displaystyle\frac{2}{(d+D)(d+4)}-\frac{5}{4(D+4)}+\frac{1}{R+4}, & R\ge 2.
\end{cases}
\]
The first expression is clearly nonnegative; and in the second, using $d\le 4$ and $R\le D$, we have
\[
\frac{2}{(d+D)(d+4)}-\frac{5}{4(D+4)}+\frac{1}{R+4}\ge \frac{1}{4(D+4)}-\frac{5}{4(D+4)}+\frac{1}{D+4}=0.
\]
So it suffices to check~\eqref{eq:one-active} when $D=4$. At $D=4$, the left-hand side of~\eqref{eq:one-active} is zero for $d=3$ (all $R$) and for $d=4$ with $R\ge 2$; indeed, asking for equality here is what explains~\eqref{eq:g3} and the second part of~\eqref{eq:g4}. For $d=4$ and $R=1$, the left-hand side equals $1/60>0$ (our choice $g(4,1)=11/60$ will be explained in Section~\ref{sec:anchors}). Finally, for $d=2$, it is zero when $R=1$, and for $2\le R\le 4$ it equals
\[
\frac{1}{12}+\frac{2-R}{(R+2)(R+4)}=\frac{(R-3)^2+23}{12(R+2)(R+4)}>0.
\]
This proves Lemma~\ref{lem:matching-cost}.
\end{proof}

\section{Edges between matching edges}\label{sec:packets}

In Lemma~\ref{lem:matching-cost}, each edge of $P$ asked for auxiliary weights from the edges joining it to other edges of $P$. Now we check that those edges can supply the amounts asked for at their ends. Because a single edge can have negative residual when both its ends are exceptional, we group these edges by pairs of $P$-edges. For disjoint vertex sets $X, Y\subseteq V_G$, let $E(X,Y)$ be the set of edges of $G$ between $X, Y$. For distinct edges $e_1, e_2\in P$, we call $E(V_{e_1},V_{e_2})$ a \emph{packet}.

\begin{lemma}\label{lem:packets}
For all distinct edges $e_1, e_2\in P$,
\begin{equation}\label{eq:packet}
E_{e_1e_2}:=\sum_{xy\in E(V_{e_1},V_{e_2})}\left(\frac{2}{d(x)+d(y)}-\lambda_x-\lambda_y\right)\ge 0.
\end{equation}
\end{lemma}

\begin{proof}
Note first that for all $d, e\in\{2,3,4\}$,
\begin{equation}\label{eq:ordinary-pair}
\frac{2}{d+e}-u_d-u_e=\frac{(4-d)(4-e)(d+e+8)}{4(d+e)(d+4)(e+4)}\ge 0.
\end{equation}
Thus every edge $xy\in E(V_{e_1},V_{e_2})$ with $\lambda_x\le u_{d(x)}$ and $\lambda_y\le u_{d(y)}$ has nonnegative residual on its own. By~\eqref{eq:u-bound} and Table~\ref{tab:weights}, this handles every edge whose ends are nonexceptional, as well as any exceptional end whose clean partner has label at least two.

The only edges left in the packet are those incident with at least one exceptional end whose clean partner has label one. Suppose first that $ac\in E(V_{e_1},V_{e_2})$ has one exceptional end $a$, with clean partner $b$ and $R_b=1$, and one nonexceptional end $c$. By Lemma~\ref{lem:signal}(a) applied from $c$ towards $b$ (since $a=\partner{b}$ and $ca\in E_G$), we have $s_c\le R_b=1$. Since $c_c>0$ (as $ca\in E_G$) and $c$ is nonexceptional, $s_c>0$; hence $s_c=1$, and so $\lambda_c=f(d(c),1)$. Put $e=d(c)\in\{2,3,4\}$; from Table~\ref{tab:weights} we get
\[
\frac{2}{3+e}-f(e,1)\ge \frac{4}{15}>\frac{3}{14}=g(3,1)
\quad \text{and} \quad
\frac{2}{4+e}-f(e,1)\ge \frac{1}{5}>\frac{11}{60}=g(4,1),
\]
so the residual of $ac$ is positive.

Now suppose that $ac\in E(V_{e_1},V_{e_2})$ has two exceptional ends, and write $e_1=ab$, $e_2=cd$. Since both $b, d$ are clean, Lemma~\ref{lem:two-to-one} implies that $bd\in E_G$ (and so $bd$ also belongs to the packet $E(V_{e_1},V_{e_2})$). We pair the residuals of $ac, bd$. Note that for all $j\in\{3,4\}$, $2\le D\le 4$, and $1\le R\le D$,
\begin{equation}\label{eq:pair-comp}
g(j,R)+f(D,R)\le u_j+u_D.
\end{equation}
Indeed, if $R\ge 2$ then~\eqref{eq:pair-comp} is immediate from $g(j,R)\le u_j$ and $f(D,R)\le u_D$. If $R=1$, then
\[
u_D-f(D,1)-\frac{3}{40}=\frac{(D-1)(D-2)}{15(D+4)}\ge 0,
\]
while $g(3,1)-u_3=3/56<3/40$ and $g(4,1)-u_4=7/120<3/40$, and~\eqref{eq:pair-comp} follows again. Applying~\eqref{eq:pair-comp} to both $ab, cd$ gives
\[
\lambda_a+\lambda_b+\lambda_c+\lambda_d\le u_{d(a)}+u_{d(b)}+u_{d(c)}+u_{d(d)};
\]
and combining this with~\eqref{eq:ordinary-pair} for $ac, bd$, we deduce that the sum of the residuals of $ac, bd$ is nonnegative.

Each edge of $P$ has at most one exceptional end, so $E(V_{e_1},V_{e_2})$ contains at most one edge with both ends exceptional. Hence at most two edges of the packet are paired in this way, and all remaining edges have nonnegative residuals. Adding over $E(V_{e_1},V_{e_2})$ gives $E_{e_1e_2}\ge 0$. This proves Lemma~\ref{lem:packets}.
\end{proof}

\section{Attachments at both ends}\label{sec:anchors}

Edges between $V_P, V_Q$ meeting $V_P$ at a nonexceptional vertex are already handled by~\eqref{eq:ordinary-bdry}. What about edges from an exceptional end to $V_Q$? For these we combine the boundary edges with their matching edge.

Let $\anch\subseteq P$ be the set of all anchored edges with one active end $a$ and one clean end $b$. By~\eqref{eq:brick-nbr}, if an exceptional end $a$ of an edge $ab\in P$ has a neighbour in $V_Q$, then $b$ also has a neighbour in $V_Q$, and so $ab\in\anch$; conversely, if $ab\in\anch$, then its active end $a$ has $c_a>0$ and $s_a=0$, and hence is exceptional. Thus, $\anch$ is just the set of edges in $P$ whose exceptional end has a neighbour in $V_Q$. For $ab\in\anch$ assigned to a brick $J\in\bricks$, put
\[
q_a=|N(a)\cap V_J| \qquad \text{and} \qquad q_b=|N(b)\cap V_J|.
\]

\begin{lemma}\label{lem:anchored}
Let $ab\in\anch$, and let $J\in\bricks$ be the brick assigned to $ab$. Then
\begin{equation}\label{eq:CabJ}
\begin{split}
C_{ab,J}:={}&\rho_{ab}+(c_a-q_a)\lambda_a+(c_b-q_b)\lambda_b\\
&+\sum_{v\in\{a,b\}}\sum_{x\in N(v)\cap V_J}\left(\frac{2}{d(v)+d(x)}-\frac{1}{r_J}\sum_{y\in N_J(x)}\frac{1}{d(x)+d(y)}\right)\ge 0.
\end{split}
\end{equation}
\end{lemma}

\begin{proof}
Let $d=d(a)$, $D=d(b)$, $R=R_b$ and $r=r_J$. As $a$ is exceptional and $b$ is clean with a neighbour in $V_J$, we have $d\in\{3,4\}$, $2\le D\le 4$ and $1\le R\le D$. By~\eqref{eq:brick-nbr}, every internal neighbour of a vertex in $N(a)\cap V_J$ lies in $N(b)\cap V_J$. In the $r$-regular graph $J$, the vertices of $N(a)\cap V_J$ have $rq_a$ internal incidences, each meeting $N(b)\cap V_J$, while $N(b)\cap V_J$ has $rq_b$ internal incidences altogether; hence $rq_a\le rq_b$ and so $q_b\ge q_a$.

In the double sum in~\eqref{eq:CabJ}, replace each $2/(d(v)+d(x))$ by $\frac{1}{r}\sum_{y\in N_J(x)}2/(d(v)+d(x))$. For every $x\in N(a)\cap V_J$ and $y\in N_J(x)$, we know that $y\in N(b)\cap V_J$, so we may pair the summand for $(x,y)$ at $v=a$ with the summand for $(y,x)$ at $v=b$. The double sum in~\eqref{eq:CabJ} then becomes
\begin{align}
&\frac{1}{r}\sum_{\substack{x\in N(a)\cap V_J\\y\in N_J(x)}}\left(\frac{2}{d+d(x)}+\frac{2}{D+d(y)}-\frac{2}{d(x)+d(y)}\right) \notag\\
&\quad +\frac{1}{r}\sum_{\substack{y\in N(b)\cap V_J\\x\in N_J(y)\setminus N(a)}}\left(\frac{2}{D+d(y)}-\frac{1}{d(x)+d(y)}\right).
\label{eq:paired-sum}
\end{align}
Because the indexing is by directed incidences of $J$, each incidence is used once, so the identity still holds even if $N(a)\cap N(b)\cap V_J\ne\emptyset$. In~\eqref{eq:paired-sum}, the paired terms have total coefficient $q_a$, and the unpaired terms have total coefficient $q_b-q_a$.

Consider first an unpaired term, for $y\in N(b)\cap V_J$ and $x\in N_J(y)\setminus N(a)$. By Lemma~\ref{lem:signal}(b) applied to $by$, we have $d(x)\ge s_b=R$, and so by~\eqref{eq:f-def} each unpaired term is at least $f(D,R)$.

Next, consider a paired term, for $x\in N(a)\cap V_J$ and $y\in N_J(x)$, and put $k=d(x)$, $l=d(y)$. Both $x, y$ have an internal neighbour in $J$ and an external neighbour ($a, b$ respectively), so $2\le k,l\le 4$; and Lemma~\ref{lem:signal}(b) applied to $by$ again gives $k\ge R$. Setting $k_0=\max(2,R)$, every paired term is at least
\[
P_1(d,D,R):=\min_{\substack{k_0\le k\le 4\\ 2\le l\le 4}}\left(\frac{2}{d+k}+\frac{2}{D+l}-\frac{2}{k+l}\right).
\]
For fixed $l$, the function $2/(d+k)-2/(k+l)$ is nondecreasing in $k$ when $l\le d$ and nonincreasing when $l\ge d$; for fixed $k\le D$, $2/(D+l)-2/(k+l)$ is nondecreasing in $l$. Since $k_0\le D$, when $d=4$ we have $l\le d$ throughout, so the minimum is attained at $(k,l)=(k_0,2)$. When $d=3$, the minimum over $l\le 3$ is likewise attained at $(k_0,2)$, while for $l=4$ it is attained at $(4,4)$. That gives
\begin{align*}
P_1(4,D,R)&=\frac{2}{4+k_0}+\frac{2}{D+2}-\frac{2}{k_0+2},\\
P_1(3,D,R)&=\min\left\{\frac{2}{3+k_0}+\frac{2}{D+2}-\frac{2}{k_0+2},\; \frac{2}{7}+\frac{2}{D+4}-\frac{1}{4}\right\}.
\end{align*}

Because $b$ is clean, $c_b=D-1$ and $\lambda_b=f(D,R)$. Adding the bound $(q_b-q_a)f(D,R)$ for the unpaired terms in~\eqref{eq:paired-sum} to $(c_b-q_b)\lambda_b$ in~\eqref{eq:CabJ}, the term $q_b$ drops out and we are left with $(D-1-q_a)f(D,R)$. Now $1\le q_a\le c_a=d-1-t_a\le 2$, so either $q_a=1$ or $q_a=2$.

First suppose that $q_a=1$. From~\eqref{eq:CabJ} and $h_a\le R$, we obtain
\[
C_{ab,J}\ge \frac{2}{d+D}+\frac{2t_a}{d+R}-1+(c_a-1)g(d,R)+(D-2)f(D,R)+P_1(d,D,R).
\]
If $d=3$, then $t_a=c_a=1$. If $d=4$, raising $t_a$ from one to two (so $c_a$ drops from two to one) adds $2/(4+R)-g(4,R)>0$ to the right-hand side. Thus it suffices to check that
\begin{align*}
A_3(D,R)&=\frac{2}{3+D}+\frac{2}{3+R}-1+(D-2)f(D,R)+P_1(3,D,R),\\
A_4(D,R)&=\frac{2}{4+D}+\frac{2}{4+R}-1+g(4,R)+(D-2)f(D,R)+P_1(4,D,R)
\end{align*}
are nonnegative for all integers $2\le D\le 4$ and $1\le R\le D$. All eighteen values are listed in Table~\ref{tab:one-attach}, and each is nonnegative. (In particular, $A_4(4,1)=0$, which explains our choice $g(4,1)=11/60$ in~\eqref{eq:g4}.)

\begin{table}[ht]
\centering
\renewcommand{\arraystretch}{1.2}
\begin{tabular}{cccc}
\toprule
$D$ & $R$ & $A_3(D,R)$ & $A_4(D,R)$ \\
\midrule
$2$ & $1$ & $113/420$ & $1/4$ \\
$2$ & $2$ & $71/420$ & $1/12$ \\
$3$ & $1$ & $1/5$ & $71/420$ \\
$3$ & $2$ & $16/105$ & $23/420$ \\
$3$ & $3$ & $11/84$ & $1/14$ \\
$4$ & $1$ & $2/105$ & $0$ \\
$4$ & $2$ & $3/35$ & $0$ \\
$4$ & $3$ & $1/10$ & $17/420$ \\
$4$ & $4$ & $3/28$ & $1/16$ \\
\bottomrule
\end{tabular}
\caption{Lower bounds for $C_{ab,J}$ when $q_a=1$.}
\label{tab:one-attach}
\end{table}

\Needspace{8\baselineskip}
Now suppose that $q_a=2$. Then $d=4$, $t_a=1$, $c_a=2$, and $D\ge 1+q_b\ge 3$. Moreover, every paired vertex $x\in N(a)\cap V_J$ and $y\in N_J(x)$ now has degree at least three: this is clear if $r\ge 2$ (as each vertex has at least two internal neighbours and one external neighbour), and if $r=1$ then $J$ is $K_2$, so $q_a=2$ implies $V_J\subseteq N(a)$ and hence $V_J\subseteq N(b)$ by~\eqref{eq:brick-nbr}, giving each vertex of $J$ its internal neighbour together with both $a, b$. Put $k_1=\max(3,R)$; we may then replace $P_1(4,D,R)$ by the minimum over $k_1\le k\le 4$ and $3\le l\le 4$, which by the same monotonicity as before is
\[
P_2(D,R)=\frac{2}{4+k_1}+\frac{2}{D+3}-\frac{2}{k_1+3}.
\]
So we need to show that
\[
B_2(D,R):=\frac{2}{4+D}+\frac{2}{4+R}-1+(D-3)f(D,R)+2P_2(D,R)\ge 0.
\]
If $D=3$, then $R\le 3$, so $k_1=3$ and $P_2(3,R)=2/7$; hence $B_2(3,R)=2/(4+R)-1/7\ge 1/7>0$. If $D=4$ and $R=1$, direct calculation gives $P_2(4,1)=5/21$ and $B_2(4,1)=53/420>0$. If $D=4$ and $R\ge 2$, then $f(4,R)\ge 1/12$, $P_2(4,R)\ge 5/21$, and $2/8+2/(4+R)-1\ge -1/2$, and so
\[
B_2(4,R)\ge -\frac{1}{2}+\frac{1}{12}+\frac{10}{21}=\frac{5}{84}>0.
\]
This proves Lemma~\ref{lem:anchored}.
\end{proof}

\section{Proof of the main theorem}\label{sec:proof}

Now we prove Theorem~\ref{thm:main}.

\begin{proof}[Proof of Theorem~\ref{thm:main}]
Let $F$ be the set of all edges $ax\in E_G$ with $a\in V_P$ and $x\in V_Q$ such that the edge of $P$ incident with $a$ is not in $\anch$. For each $ax\in F$, the vertex $a$ is nonexceptional by~\eqref{eq:brick-nbr} and the definition of $\anch$, so~\eqref{eq:ordinary-bdry} holds for $ax$. Let $\Xi$ be the sum of: $S_{ab}$ over all $ab\in P\setminus\anch$; $C_{ab,J}$ over all $ab\in\anch$ (where $J$ is the brick assigned to $ab$); $E_{e_1e_2}$ over all unordered pairs of distinct edges $e_1,e_2\in P$; the left-hand side of~\eqref{eq:brick-balance} over all bricks $J\in\bricks$; and $2/(d(a)+d(x))-\lambda_a-\phi_{J,x}$ over all $ax\in F$ (where $J\in\bricks$ is the brick containing $x$). By Lemmas~\ref{lem:matching-cost}, \ref{lem:packets} and~\ref{lem:anchored}, together with~\eqref{eq:brick-balance} and~\eqref{eq:ordinary-bdry}, every term in this sum is nonnegative, and so $\Xi\ge 0$.

We claim that all terms in $\lambda_a$ ($a\in V_P$) and $\phi_{J,x}$ ($J\in\bricks$, $x\in V_J$) cancel in $\Xi$. To see this, let $a\in V_P$ and let $b=\partner{a}$. If $ab\notin\anch$, then $c_a\lambda_a$ comes from $S_{ab}$, while every edge from $a$ to $V_M\setminus\{a,b\}$ contributes $-\lambda_a$ once (in a packet if the other end lies in $V_P$, and in the sum over $F$ if it lies in $V_Q$); as there are $c_a$ such edges, $\lambda_a$ cancels. If $ab\in\anch$ is assigned to a brick $J$, then no edge at $a$ is in $F$, and all neighbours of $a$ in $V_Q$ lie in $V_J$ by Lemma~\ref{lem:same-brick}. Hence $a$ has $c_a-q_a$ edges to $V_P\setminus\{a,b\}$, and their $-\lambda_a$ terms in the packets cancel $(c_a-q_a)\lambda_a$ in $C_{ab,J}$. For a brick $J\in\bricks$ and $x\in V_J$, the term $(d(x)-r_J)\phi_{J,x}$ in~\eqref{eq:brick-balance} likewise cancels the $-\phi_{J,x}$ terms from the $d(x)-r_J$ edges joining $x$ to $V_P$: each such edge $ax$ ($a\in V_P$) gives $-\phi_{J,x}$ in~\eqref{eq:CabJ} if $a\partner{a}\in\anch$ (where the inner sum is $\phi_{J,x}$) and in the sum over $F$ otherwise.

It remains to collect the non-auxiliary terms of $\Xi$. The $-1$ terms in $\rho_{ab}$ ($ab\in P$) and the $-|V_J|/2$ terms in~\eqref{eq:brick-balance} ($J\in\bricks$) sum to $-|P|-|Q|=-|M|$. Because $U_M$ is independent and no vertex of $V_Q$ is adjacent to $U_M$, every edge of $G$ either has both ends in $V_M$ or goes between $V_P$ and $U_M$. Take an edge $uv$ with $u,v\in V_M$; we check that $2/(d(u)+d(v))$ occurs once in $\Xi$. If $u,v\in V_P$, then it shows up in $\rho_{uv}$ when $uv\in P$ and in a packet otherwise. If $u,v\in V_Q$, then $u,v$ lie in the same brick $J$ and $uv$ is counted in~\eqref{eq:brick-balance}. If one end is in $V_P$ and the other is in $V_Q$, say $a\in V_P$ and $x\in V_Q$, then $ax$ is counted in~\eqref{eq:CabJ} when $a\partner{a}\in\anch$ (as $x$ lies in the brick assigned to $a\partner{a}$ by Lemma~\ref{lem:same-brick}) and in the sum over $F$ otherwise. An edge $vz$ with $v\in V_P$ and $z\in U_M$ is also counted once, in $\rho_{v\partner{v}}$, except that $2/(d(v)+h_v)$ appears there in place of $2/(d(v)+d(z))$. Adding back the missing difference over all such edges, we get
\begin{equation}\label{eq:final}
H(G)-|M|=\Xi+\sum_{v\in V_P}\sum_{z\in N(v)\cap U_M}\left(\frac{2}{d(v)+d(z)}-\frac{2}{d(v)+h_v}\right).
\end{equation}
Both terms on the right-hand side of~\eqref{eq:final} are nonnegative, as $\Xi\ge 0$ and $d(z)\le h_v$ for all $v\in V_P$ and $z\in N(v)\cap U_M$. (Any isolated vertex of $G$ sits in $U_M$ and contributes nothing to either side.) Thus $H(G)\ge |M|=\mustar(G)$. This proves Theorem~\ref{thm:main}.
\end{proof}

\end{document}